\documentclass[10pt,a4paper]{article}
\usepackage[latin1]{inputenc}
\usepackage{amsmath}
\usepackage{amsfonts}
\usepackage{amssymb}

\title { In\'egalit\'es de Harnack et ph\'enom\`ene de concentration.}

\author{Samy Skander Bahoura \footnote { Adresses e-mails: samybahoura@yahoo.fr, bahoura@ccr.jussieu.fr }}

\date{ \small { \it { Universit\'e de Patras, D\'epartement de Math\'ematiques, Patras 26500 Gr\`ece }}}

\begin{document}

\maketitle

\bigskip

 \hrule 

\bigskip

{\bf R\'esum\'e}. Nous prouvons des in\'egalit\'es concernant le produit $ \sup u \times \inf u $ pour des op\'erateurs elliptiques d'ordre 2 et 4. Ces in\'egalit\'es et le ph\'enom\`ene de concentration nous permettent d'obtenir le comportement asymptotique des solutions de ces EDP.

\bigskip

{\bf Abstract}. We proove some inequalities concerning the product $ \sup u \times \inf u \leq c $, $ \sup u \times \inf u \geq c $  for some elliptic operators of order 2 and 4. Those inequalities and the concentration phenomena we can describe the asymptotic behavior of those PDE solutions.

\bigskip

{\it Mots-cles}: $ \sup \times \inf $, laplacien et bilaplacien, point de concentration, comportement asymptotique. 


\bigskip

\hrule

\bigskip

Cet article correspond \`a la Note aux Comptes Rendus Math\'ematiques de l'Academie des Sciences [B 3].

\bigskip

Dans la suite nous notons le laplacien g\'eom\'etrique par $ \Delta=-\nabla^i\nabla_i $.

\bigskip

On s'occupe de certaines in\'egalit\'es de Harnack de type $ \sup \times \inf $ et leurs applications aux ph\'enom\`enes de concentration dans le cas d'op\'erateurs elliptiques d'ordre 2. Ce type de probl\`emes est bien connu voir par exemple, [B1], [B2] [C-L], [DLN], [GT], [H], [He] [L1], [L2], [M] et les r\'esultats obtenus utilisent les t\'echniques de sym\'etrie voir [GNN].

\bigskip

Les ph\'enom\`enes de concentration et leurs cons\'equences ont beaucoup \'et\'e \'etudi\'ees, precis\'ement dans la recherche de meilleurs constantes dans les in\'egalit\'es de Sobolev voir par exemple [Au], [Au, D, He], [DHR] et [He,V]. En ce qui nous concerne, ce type d'in\'egalit\'es ( $ \sup \times \inf $ ) et le ph\'enom\`ene de concentration, nous permettent de d\'ecrire le comportement asymptotic de certaines solutions d'EDP.

\bigskip

Consid\'erons une suite de r\'eels positifs $ (\epsilon_i)_{i\geq 0} $ avec $ \epsilon_i \to 0 $ et une suite de fonctions $ v_{\epsilon_i} >0 $ sur $ {\mathbb S}_n $, telles que:

$$ \Delta v_{\epsilon_i}+\dfrac{n(n-2)}{4}v_{\epsilon_i}=\dfrac{n-2}{4(n-1)}V_{\epsilon_i}{v_{\epsilon_i}}^{N-1-\epsilon_i} , $$

avec, $ N=\dfrac{2n}{n-2} $, $ 0 < a \leq V_{\epsilon_i}(x)\leq b <+\infty, \forall \,\, x\in {\mathbb S}_n $ et $ ||\nabla V_{\epsilon_i}||_{\infty}\leq A $.

\bigskip

On suppose que pour tout $ i $ que les fonctions $ V_{\epsilon_i} $ et $ v_{\epsilon_i} $ sont r\'eguli\`eres.

\bigskip

{\it Th\'eor\`eme1. Sous ces hypoth\`eses la suite $ (v_{\epsilon_i}) $ v\'erifie,

$$ {\epsilon_i}^{(n-2)/2} \left ( \sup_{{\mathbb S}_n} v_{\epsilon_i} \right )^{1/4} \times \inf_{{\mathbb S}_n} v_{\epsilon_i} \to 0.\qquad (1) $$}

On s'int\'eresse aussi au probl\`eme suivant:

$$ \Delta u_i=V_i(|x|){u_i}^{N-1-\epsilon_i} \,\,\, {\rm dans }\,\, B_1(0) \,\,\, {\rm et} \,\,\, u_i=0 \,\, {\rm sur} \,\, \partial B_1(0). $$

O\`u $ B_1(0) $ est la boule unit\'e de $ {\mathbb R}^n $ et $ V_i $ est une fonction d\'ecroissante de $ |x| $ qui v\'erifie $ 0 < a\leq V_i(|x|)\leq b <+\infty $ pour tout $ x $. On suppose que les fonctions $ u_i $ et $ V_i $ sont r\'eguli\`eres pour tout $ i $.

\bigskip

{\it Th\'or\`eme 2. Pour tout compact $ K $ de $ B_1(0) $,

$$ \sup_{B_1(0)} u_i\times \inf_K u_i \geq c=c(a,b,K,\Omega,n) >0 . $$}

{\underbar {\bf Application:}}

\bigskip

Consid\'erons le probl\`eme suivant:

$$ \Delta u_{\epsilon_i}=n(n-2){u_{\epsilon_i}}^{N-1-\epsilon_i},\,\,\, u_{\epsilon_i}>0 \,\,\, {\rm dans}\,\,\, \Omega \,\, {\rm et} \,\, u_{\epsilon_i}=0 \,\, {\rm sur} \,\, \partial \Omega \qquad (E), $$

avec, $ \Omega $ un ouvert born\'e de $ {\mathbb R}^n $.

\bigskip

Pour ce type d'\'equation il existe de nombreux r\'esultats de compacit\'e et de comportement asymptotique voir par exemple, [BP], [H], [He]

\bigskip

Nous avons,

\bigskip

{\it Th\'eor\`eme 3. 1) Il existe $ c_1=c_1(n, \Omega)>0, c_2=c_2(n, \Omega)>0 $ telles que:

$$ c_2 \leq ||u_{\epsilon_i}||_{H^1(\Omega)} \leq c_1. $$

2) Si $ \Omega $ est \'etoil\'e, alors, il existe une sous-suite $ (u_{\epsilon_j}) $ pour laquelle, il existe un $ m \in {\mathbb N}^* $ et un nombre fini de points de concentrations $ x_1,x_2,...,x_m\in \Omega $ tels que:

$$ i) \lim_{\epsilon_j \to 0} u_{\epsilon_j} = 0 \,\,\, {\rm dans } \,\, {\cal C}_{loc}^2(\bar \Omega-\{x_1,\ldots x_m\}), $$

$ \forall \,\, k\in \{1,\ldots, m\}, \,\, \exists \,\, (x_{j,k}) $ avec , $ x_{j,k} \to x_k $ et $ u_{\epsilon_j}(x_{j,k})\to+\infty $.

$$ ii) \lim_{\epsilon_j \to 0}  u_{\epsilon_j}^{N-\epsilon_j} =\sum_{i=1}^m \mu_i\delta_{x_i} \,\,\, {\rm avec } \,\, \mu_i\geq \dfrac{\omega_n}{2^n}. $$

Ici la convergence est au sens des distributions.

\bigskip

 iii)  Pour tout compact $ K $ de $ \Omega-\{x_1,\ldots,x_m \} $, il existe une constante positive $ c=c(K,\Omega,n)>0 $ telle que:

$$ \sup_{\Omega} u_{\epsilon_j} \times \sup_K u_{\epsilon_j} \leq c. $$

 iv)  Il existe un voisinage $ \omega $ du bord  $ \partial \Omega $ et une constante positive $ \bar c=\bar c(\omega,\Omega,n) $ tels que:

$$ \sup_{\Omega} u_{\epsilon_j} \times \sup_{\omega} u_{\epsilon_j} \leq \bar c. $$
 
 v)  il existe deux constantes positives, $ \beta_1 $ et $ \beta_2 $ telles que:

$$ \beta_1 \leq \epsilon_j \left ( \sup_{\Omega} u_{\epsilon_j} \right )^2 \leq \beta_2, $$

plus pr\'ecis\'ement, il existe une fonction $ g \in {\cal C}^2(\partial \Omega) $, telle que,

$$  \epsilon_j \left ( \sup_{\Omega} u_{\epsilon_j} \right )^2 \to \dfrac{ c_n \int_{\partial \Omega} <x|\nu(x)>[ \partial_{\nu} g(\sigma)]^2 d\sigma }{\sum_{k=1}^m \mu_k}. $$

  vi) il existe $ m $ r\'eels positifs $ \gamma_1, \ldots, \gamma_m $, $ \gamma_k\geq n(n-2)\dfrac{\omega_n}{2^n} $, $ k\in\{1,\ldots,m \} $, tels que:
  
  $$ \sup_{\Omega} u_{\epsilon_j} \times u_{\epsilon_j}(x) \to \sum_{k=1}^m \gamma_k G(x_k,x) \,\, {\rm dans } \,\, {\cal C}_{loc}^2(\bar \Omega-\{x_1,\ldots,x_m\}), $$
 
  o\`u $ G $ est la fonction de Green du laplacien avec condition de Dirichlet. On peut prendre, $ g=\sum_{k=1}^m \gamma_k G(x_k,.) $ dans le v).}
 
\bigskip

{\bf Remarque:} On verra que dans le cas o\`u $ \Omega $ n'est pas \'etoil\'e, les points i), ii), iii), iv) et vi) se conservent.

\bigskip

Concernant certains op\'erateurs d'ordre 4, il existe des r\'esultats comme [C 1], [ C 2], [C-G] et [V]. Nous nous occupons maintenant de savoir si pour certains d'entre eux des minorations du produit $ \sup \times \inf $ sont possibles.

\bigskip

Sur une vari\'et\'e Riemannienne compacte $ (M,g) $ de dimension $ n\geq 5 $, on consid\`ere l'\'equation suivante:

$$ \Delta^2 u_i+b\Delta u_i+cu_i=V_i{u_i}^{(n+4)/(n-4)},\,\, u_i > 0 \,\, {\rm sur} \,\, M  \qquad \qquad (E')$$

avec, $ b,c>0 $, $ c\leq \dfrac{b^2}{4} $ et $ 0 \leq V_i(x)\leq A $.

La condition $ 0 <c \leq \dfrac{b^2}{4} $ est tr\'es utile pour obtenir notre estimation, elle permet d'avoir une fonction de Green avec d'int\'er\'essantes propri\'et\'es et est utilis\'ee pour appliquer le principe du maximum voir [C 2].

\bigskip

{\it Th\'eor\`eme 4. Il existe une constante positive $ k=k(b,c,A,M,g)  $, telle que,

$$ \sup_M u_i \times \inf_M u_i \geq k\,\,\forall \,\, i, $$

o\`u $ u_i $ est solution de $ (E') $ }

\bigskip

Autre r\'esultat sur un ouvert $ \Omega $ strictement convexe de $ {\mathbb R}^n $ avec $ n\geq 5 $, consid\'erons l'\'equation:

$$ \Delta^2 u_{\epsilon}={u_{\epsilon}}^{p-\epsilon},\,\,\,  u_{\epsilon} >0 \,\,\, {\rm dans} \,\,\Omega,\,\,\, {\rm et}\,\, u_{\epsilon}=\Delta u_{\epsilon}=0\,\,\,{\rm sur} \,\,\partial \Omega, \qquad (E'') $$

avec, $ p=\dfrac{n+4}{n-4} $ et $ 0 < \epsilon \leq \dfrac
{4}{n-4} $.

\bigskip

%{\it Th\'eor\`eme 5. Il existe $ \delta=\delta(n,\Omega)>0 $ tel que pour tout solution $ u_{\epsilon} $ de $ (E'') $, on ait:

%$$ d(x_{\epsilon},\partial \Omega)\geq \delta,\,\,\, {\rm avec }\,\,\, u_{\epsilon}(x_{\epsilon})=\max_{\Omega} u_{\epsilon} . $$}

{\it Th\'eor\`eme 5. Pour tout compact $ K\subset \Omega $, il existe $ c=c(K,\Omega,n)>0 $ telle que pour toute solution $ u_{\epsilon} $ de $ (E'') $, on ait:

$$ \sup_{\Omega} u_{\epsilon} \times \inf_{K} u_{\epsilon} \geq c . $$ }

Pour les op\'erateur d'ordre 4, il existe une identit\'e de Pohozaev comme pour ceux d'ordre 2( [P]), voir [C-G], et l'\'equation $ (E'') $, avec $ \epsilon=0 $, ne poss\`ede pas de solutions lorsque l'ouvert $ \Omega $ est \'etoil\'e. Ceci nous pousse \`a \'etudier $ (E'') $ avec $ \epsilon >0 $.

\smallskip

Pour l'\'equation $ (E'') $, il est important de remarquer que les points $ x_{\epsilon} $ o\`u les  solutions $ u_{\epsilon} $ sont maximum, restent loin du bord $ \partial \Omega $. Ceci est important pour estimer ces solutions pr\'es du bord. Ce fait important est r\'ealis\'e pour les ouverts $ \Omega $ strictement convexes, la m\'ethode utilis\'ee est celle de sym\'etrie de [GNN] et la transformation de Kelvin. Pour le laplacien d'ordre 2, l'\'equation reste invariante par la transformation de Kelvin, ce qui n'est pas le cas pour $ (E'') $  et la condition de stricte convexit\'e du domaine peut attenuer la non invariance par cette transformation, voir [C-G].

\newpage

\underbar{\bf {Preuve du Th\'eor\`eme 1:}}

\bigskip

Supposons par l'absurde que $ (1) $ ne soit pas vraie, alors:

$$ \limsup_{\epsilon_i \to 0} \left [ {\epsilon_i}^{(n-2)/2}  \left ( \sup_{{\mathbb S}_n} v_{\epsilon_i} \right )^{1/4} \times \inf_{{\mathbb S}_n} v_{\epsilon_i} \right ] \geq c >0. $$

Soit $ x_{\epsilon_i} $ le point o\`u $ v_{\epsilon_i} $ atteint son maximum. On consid\`ere la projection st\'eographique de p\^ole $ y_{\epsilon_i} $ point diam\`etralement oppos\'e \`a $ x_{\epsilon_i} $.

\bigskip

On a:

\smallskip

Si $ x=(x_1, \ldots, x_n,x_{n+1}) \in {\mathbb S}_n $ et $ y=(y_1,\ldots, y_n)\in {\mathbb R}^n $ l'image de $ x $ par la projection st\'er\'eographique, alors: 

$$ x_i=\dfrac{2y_i}{1+|y|^2}, \,\,\, 1 \leq i \leq n, x_{n+1}=\dfrac{|y|^2-1}{|y|^2+1} , $$

et,

$$ y_i=\dfrac{x_i}{1-x_{n+1}},\,\,\, 1 \leq i \leq n .$$

$ g_0=\left (\dfrac{2}{1+|y|^2} \right )^2 dy^2,\,\,\, {\rm et} \,\,\, u_{\epsilon_i}=\left (\dfrac{2}{1+|y|^2} \right )^{(n-2)/2} v_{\epsilon_i} .$

\smallskip

Soit $ {\cal E} $ la m\'etrique euclidienne et $ \Delta_{\cal E} $ le laplacien pour cette m\'etrique, on obtient:

$$ \Delta_{\cal E} u_{\epsilon_i}=\dfrac{n-2}{4(n-1)}H^{\epsilon_i} V_i{u_{\epsilon_i}}^{N-1-\epsilon_i},\,\,\, {\rm  avec } \,\,\, H(y)=\left (\dfrac{2}{1+|y|^2} \right )^{(n-2)/2}. $$

et,

$$ u_{\epsilon_i}(0)=2^{(n-2)/2}\max_{{\mathbb S}_n} v_{\epsilon_i}=\max_{B_1(0)} u_{\epsilon_i},\,\,\, {\rm et } \,\,\, \inf_{B_1(0)} u_{\epsilon_i} \geq \inf_{{\mathbb S}_n} v_{\epsilon_i}, $$

On en d\'eduit que,

$$ \,\,\,\limsup_{\epsilon_i \to 0} \left [ {\epsilon_i}^{(n-2)/2} \left [ u_{\epsilon_i}(0)\right ]^{1/4} \times \inf_{B_1(0)} u_{\epsilon_i} \right ] \geq \tilde c >0. $$

On se retrouve dans le cas du Th\'eor\`eme 1 de  [B 1], les \'etapes de la preuve de ce Th\'eor\`eme 1 se conservent ici, mais il faut v\'erifier si le lemme 2 de l'\'etape 2-3 se conserve.

\bigskip

On a, $ a_{\epsilon_i}=0 $ pour tout entier $ i $ et,

$$ \tilde V_{\epsilon_i}(t,\theta)=e^{(n-2)\epsilon_i t/2}[H(t)]^{\epsilon_i}V_{\epsilon_i}(e^t\theta), \,\,\, {\rm avec }\,\,\, H(t)=\left ( \dfrac{2}{1+e^{2t}} \right )^{(n-2)/2}. $$

Comme,

$$ {\epsilon_i}^{[2/(n-2)-\epsilon_i/2]^{-1}}\left [u_{\epsilon_i}(0)\right ]^{1/4}\geq \tilde c>0 , $$

ce qui s'\'ecrit,

$$ \log \epsilon_i\geq -\dfrac{1}{4}\left ( \dfrac{2}{n-2}-\dfrac{\epsilon_i}{2} \right )\log u_{\epsilon_i}(0)+\log {\tilde {\tilde c}}=t_i . $$

et donc, 

$$ \partial_t \tilde V_{\epsilon_i} \geq 0,\,\,\, {\rm sur} \,\,\, ]-\infty, t_i]\times {\mathbb S}_{n-1}. $$

Donc, le lemme 2 de l'\'etape 2-3 du Th\'eor\`eme 1 de [B 1] se conserve et la conclusion de la preuve de ce Th\'eor\`eme 1 reste la m\^eme ici, \`a savoir,

$$ [u_{\epsilon_i}(0)]^{1/4}\inf_{B_1(0)} u_{\epsilon_i} \leq \bar c ,\,\,\, {\rm pour \,\, tout } \, \, i , $$

Or, d'apr\`es notre hypoth\`ese de d\'epart,

$$ [u_{\epsilon_i}(0)]^{1/4} \inf_{B_1(0)} u_{\epsilon_i} \geq \dfrac{\tilde c}{{\epsilon_i}^{(n-2)/2}} \to +\infty .$$

\bigskip

{\underbar {\bf Preuve du Th\'eor\`eme 2:}}

\smallskip

D'apr\`es le Th\'eor\`eme 1' de Gidas-Ni-Nirenberg [GNN], la m\'ethode moving-plane assure que $ u_i $ est radiale. Consid\'erons la fonction de Green $ G $ du laplacien, on peut ecrire:

$$ \max_{B_1(0)} u_i=u_i(0)=\int_{B_1(0)} G(0,y)V_i(|x|){u_i(y)}^{N-1-\epsilon_i} \leq b\left [\int_{B_1(0)} G(0,y)dy \right ](\max_{B_1(0)} u_i)^{N-1-\epsilon_i} ,$$

Ce qui donne:

$$ [u_i(0)]^{4/(n-2)-\epsilon_i} \geq c>0, $$

avec, $ c=\dfrac{1}{b\left [\int_{B_1(0)} G(0,y)dy \right ]} $.

\smallskip

Comme dans le Th\'eor\`eme 2 de [B-2], on utilise la représentation par la fonction de Green pour prouver que l'\'energie tend vers 0, puis par le proc\'ed\'e d'it\'eration de Moser, on prouve que $ u_i(0) \to 0 $ ce qui est contradictoire.

\bigskip

{\underbar {\bf Preuve du Th\'eor\`eme 3:}}

\smallskip

%On sait qu'il existe $ \delta=\delta(n,\Omega)>0 $ tel que:

%$$ d(x_{\epsilon},\partial \Omega)\geq \delta, $$

%avec $ u_{\epsilon}(x_{\epsilon})=\max_{\Omega} u_{\epsilon} $.

%Soit $ G $ la fonction de Green du laplacien avec condition de Dirichelet, alors:

%$$ u_{\epsilon}(x)=\int_{\Omega} G(x,y) u_{\epsilon}^{N-1-\epsilon}(y) dy. $$

\bigskip

\underbar {\bf Preuve de 1):}

D'apr\`es l'in\'egalit\'e de Harnack du type $ \sup \times \inf $ (voir [C-L], [L1] ou adapter la m\'ethode qui se trouve dans  [B 1] avec les $ V_i \equiv 1 $ ), pour tout compact $ K $ d'un ouvert $ \Omega_0 \subset \subset \Omega $, il existe une constante positive  $ \bar c= \bar c(K, \Omega_0, \Omega, n) $ telle que :

$$ \sup_K u_{\epsilon} \times \inf_{\Omega_0} u_{\epsilon} \leq c , $$

%$$ u_{\epsilon}(x_{\epsilon})u_{\epsilon}(y_{\epsilon})=\sup_{ \{ x,d(x,\partial \Omega) \leq \delta /2\} } u_{\epsilon} \times \inf_{ \{x, d(x,\partial \Omega)\leq \delta \}} u_{\epsilon} \leq c(\delta,\Omega,n) . $$

On sait que (voir [H]  pages 164-165) pour $ \delta=\delta_0 >0 $ assez petit:

$$ u_{\epsilon}(x) \leq M=M(\delta,\Omega,n) \,\,\, \forall \,\, x\in \{y, d(y, \partial \Omega)\leq \delta \}. $$

$ \delta_0 $ d\'epend de $ \Omega $ mais pas de $ \epsilon >0 $. 

\bigskip

On prend, $ \Omega_0=\{x, d(x,\partial \Omega)> \dfrac{\delta_0}{2} \} $ et $ K=\{x, d(x,\partial \Omega)\geq \dfrac{2\delta_0}{3} \} $.

\bigskip

Soit $ G $ la fonction de Green du laplacien avec condition de Dirichelet, alors:

$$ u_{\epsilon}(x)=\int_{\Omega} G(x,y) u_{\epsilon}^{N-1-\epsilon}(y) dy. $$

%Donc,

%$$ c(\delta, \Omega, n) \geq  \int_{\Omega} G(y_{\epsilon},y)u_{\epsilon}^{N-\epsilon}(y)dy . $$

En appliquant le principe du maximum \`a $ G $, on a pour $ \delta_0 >0 $,

$$ G(x,y)\geq c_1=c_1(\delta_0,\Omega,n)>0\,\,\, {\rm sur} \,\, \{x, d(x,\partial \Omega)\geq 2\delta_0/3 \} \times \{y, d(y,\partial \Omega)\geq \delta_0/2 \}  , $$

ce qui donne,

$$ \bar c \geq \sup_K u_{\epsilon} \times \inf_{\Omega_0} u_{\epsilon} \geq \int_{\{x, d(x,\partial \Omega)\geq 2\delta_0/3 \} } u_{\epsilon}^{N-\epsilon}(y)dy . $$

%On sait que (voir [H]) pour $ \delta=\delta_0 >0 $ assez petit:

%$$ u_{\epsilon}(x) \leq M=M(\delta,\Omega,n) \,\,\, \forall \,\, x\in \{y, d(y, \partial \Omega)\leq \delta \}. $$

%$ \delta_0 $ d\'epend de $ \Omega $ mais pas de $ \epsilon >0 $. 

\smallskip

Comme,

$$ ||u_{\epsilon}||_{N-\epsilon}^{N-\epsilon}=\int_{\{x, d(x,\partial \Omega)\leq 2\delta_0/3 \}} u_{\epsilon}^{N-\epsilon}(y)dy+\int_{\{x, d(x,\partial \Omega)\geq 2\delta_0/3 \} } u_{\epsilon}^{N-\epsilon}(y)dy. $$

On en d\'eduit que:

$$ ||u_{\epsilon}||_{N-\epsilon}^{N-\epsilon} \leq c(n,\Omega). $$

En multipliant $ (E) $ par $ u_{\epsilon} \in H_0^1(\Omega) $ et en int\`egrant par parties, on obtient:

$$ ||u_{\epsilon}||_{H_0^1(\Omega)}\leq c(n,\Omega) . $$

D'autre part, en utilisant l'injection de Sobolev, puis en multipliant l'\'equation $ (E) $ par $ u_{\epsilon} $, en int\'egrant par parties et en utilisant l'in\'egalit\'e de Holder, on obtient:

$$ K_1||u_{\epsilon}||_{N-\epsilon}^2\leq K_2||u_{\epsilon}||_N^2\leq \int_{\Omega} |\nabla u_{\epsilon}|^2=n(n-2)\int_{\Omega} u_{\epsilon}^{N-\epsilon}=n(n-2)||u_{\epsilon}||_{N-\epsilon}^{N-\epsilon}, $$

Ce qui donne,

$$ ||u_{\epsilon}||_{N-\epsilon} \geq \tilde K_1>0 \,\,\, {\rm et} \,\, ||u_{\epsilon}||_{H^1} \geq \tilde K_2>0 . $$

\underbar {\bf Preuve de 2):}

\bigskip

\underbar {\bf Point i):}

Comme $ u_{\epsilon} $ est born\'ee dans $ H^1 $, on peut en extraire une sous-suite $ u_{\epsilon'} \to u $ et la convergence est presque partout, dans $ L^{N-(1/2)} $ fort et dans $ H^1 $ faible.

\bigskip

La fonction $ u\geq 0 $ v\'erifie, $ \Delta u=n(n-2)u^{N-1} $ et $ u=0 $ sur  $ \partial \Omega $. De plus $ ||u||_{H^1(\Omega)} \leq \liminf_{\epsilon \to 0} ||u_{\epsilon_i}||_{H^1(\Omega)} \leq c(n,\Omega) $.

\bigskip

{\underbar {\bf 1') R\'egularit\'e de u au bord:}}

\bigskip

On sait (voir [H] pages 164-165), qu'il existe un voisinage $ \omega $ du bord $\partial \Omega $ et constante positive $ M $ telle que:

$$ u_{\epsilon_i}(x)\leq M \,\, \forall \,\, x \in \omega. $$

En utilisant le lemme 2 de [H] (ou le lemme 8 de [He]), on en d\'eduit que la suite $ (u_{\epsilon_i}) $ est uniform\'ement born\'ee dans $ {\cal C}^{1,\alpha} $ au voisinage du bord. En utilisant le th\'eor\`eme d'Ascoli, la suite $ (u_{\epsilon_i}) $ converge vers une fonction $ {\cal C}^1 $ et donc $ u $ est $ {\cal C}^1 $ au voisinage du bord.

\bigskip

{\underbar {\bf 2') R\'egularit\'e de u \`a l'int\'erieur de $ \Omega $:}}

\bigskip

En utilisant les arguments de r\'egularit\'e locale de Trudinger (voir [T] th\'eor\`eme 3 et lemme page 268) et le th\'eor\`eme de Ladyzenskaya-Ural'ceva (voir [Au] th\'eor\`eme 4.40), on en d\'eduit que $ u $ est au moins de classe $ {\cal C}^2(\Omega) $.

\bigskip

Finalement $ u\in {\cal C}^2(\bar \Omega) $. D'apr\'es le principe du maximum $ u\equiv 0 $ ou $ u>0 $. Comme $ \Omega $ est \'etoil\'e, en utilisant la formule de Pohozaev [P], on conclut que $ u\equiv 0 $.

\bigskip

\underbar {\bf Remarque:} On peut se passer du fait que $ \Omega $ est \'etoil\'e. On sait que les points $ x_{\epsilon_i} $ o\`u $ u_{\epsilon_i} $ est maximum restent loin du bord, alors, soit $ u_{\epsilon_i} $ est uniform\'ement born\'e dans $ {\cal C}^0(\bar \Omega) $ et on a donc un r\'esultat de compacit\'e, soit, il existe une sous-suite de $ u_{\epsilon_i} $ not\'ee encore $ u_{\epsilon_i} $ telle que $ u_{\epsilon_i}(x_{\epsilon_i}) \to +\infty $ et dans ce cas l'in\'egalit\'e de Harnack du type $ \sup \times \inf $ et les estimations de Schauder ( au voisinage du bord), nous permettent d'avoir l'existence d'un point $ a_0 $ tr\'es proche du bord sans \^etre sur le bord ( o\`u $ u_{\epsilon_i} $ est uniform\'ement born\'ee) tel que $ u(a_0)=0 $, le principe du maximum appliqu\'e \`a $ u $ implique que $ u\equiv 0 $.

\bigskip

On d\'efinit un point de concentration $ x_0 $ comme suit,

$$ \forall \,\, \delta >0 \,\, \liminf_{\epsilon_i \to 0} \int_{B(x_0,\delta)\cap \Omega} u_{\epsilon_i}^{N-\epsilon_i}(x)dx >0 . $$

On prend la $ \liminf $ au lieu de $ \limsup $, dans la d\'efinition d'un point de concentration, pour pouvoir compter tous les points de concentrations, on verra plustard qu'il est plus simple de compter les points de concentrations ( qui serons en nombre fini).

\bigskip

En utilisant le Lemme 1 dans [He] (avec la modification $ \limsup $ devient $ \liminf $), on en d\'eduit que si $ x_0 $ est un point de concentration, alors:

$$ \liminf_{\epsilon \to 0} \int_{B(x_0,\delta)\cap \Omega} u_{\epsilon}^{N-\epsilon}(y)dy \geq \dfrac{\omega_n}{2^n}. $$

\bigskip

Comme $ u_{\epsilon} (x)\leq M $ sur $ \{x,d(x,\partial \Omega) \leq \delta \} $, par les estimations elliptiques on en d\'eduit que $ ||\nabla u_{\epsilon}||_{L^{\infty}(\omega)} \leq M'$. o\`u $ \omega $ est un voisisnage du bord $ \partial \Omega $. On conclut gr\^ace au Th\'eor\`eme d'Ascoli que $ u_{\epsilon} $ converge uniform\'ement vers $ 0 $ sur un voisinage du bord et donc il n'y a pas de points de concentration au voisinage du bord. De plus en utilisant le le lemme 2 de [H] ou lemme 8 de [He] et le th\'eor\`eme 6.6 de [GT], on d\'eduit que $ || u_{\epsilon}||_{{\cal C}^2(\omega)} $ converge vers 0 au voisinage $ \omega $ du bord.( on prend une couronne $ C $ voisinage du bord, les bords $ \partial C $ de $ C $ sont, $ \partial \Omega $ et $ \partial (\Omega_{\epsilon_0}) $, on prend $ \eta $ une fonction telle que $ \eta \equiv 0 $ sur $ \partial (\Omega_{\epsilon_0}) $ et $ \eta \equiv 1 $ dans un voisinage de $ \partial \Omega $, puis on s'occupe de l'\'equation $ \Delta (\eta u_{\epsilon})=f_{\epsilon} \in {\cal C}^{0,\beta}(C) $ et $ \eta u_{\epsilon}=0 $ sur $ \partial C $ ).

\bigskip

Comme $ (u_{\epsilon_i}) $ est born\'ee dans $ H^1(\Omega) $, par le th\'eor\`eme de Rellich-Kondrachov, pour tout $ q\in [1,N[ $, il existe une sous- suite not\'ee encore $ (u_{\epsilon_i}) $ qui converge vers $ 0 $ dans $ L^q $. En prenant une suite $ q \to N $ et en utilisant le proc\'ed\'e diagonal, on d\'eduit qu'il existe une sous-suite not\'ee encore $ (u_{\epsilon_i})$ qui conevrge vers $ 0 $ dans tous les $ L^q $ avec $ q\in [1, N[ $.

\bigskip

En premi\`ere conclusion, on obtient, une sous-suite not\'ee $ ( u_{\epsilon_i}) $ telle que $ ||u_{\epsilon_i}||_{{\cal C}^2(\omega)} \to 0 $ et $ \forall \,\, q\in [1,N[ $, $ ||u_{\epsilon_i}||_{L^q(\Omega)} \to 0 $, avec $ \omega $ un voisinage du bord.

\bigskip

Comme $ u_{\epsilon_i} $ est born\'ee dans $ L^N $, on en d\'eduit que si il y a des points de concentration alors ils sont en nombre fini $ x_1,\ldots x_m $. De plus $ m $ est born\'e par un r\'eel ne d\'ependant que de $ \Omega $ et $ n $.

\bigskip

On va voir que la suite $ u_{\epsilon_i} $ a au moins un point de concentration. Supposons le contraire, soit alors $ K $ un compact de $ \Omega $ tel que $ \Omega-K\subset \omega $, o\`u $ \omega $ est un voisinage du bord tel que, $ ||u_{\epsilon_i}||_{{\cal C}^2(\omega)} \to 0 $, alors $ \sup_K u_{\epsilon_i}=u_{\epsilon_i}(y_{\epsilon_i}) $ avec $ y_{\epsilon_i} \to y\in K $. Apr\'es passage \`a une sous-suite, il existe  $ \delta_y>0 $ tel que, $ \liminf_{\epsilon_i \to 0} \int_{B(y,\delta_y)} u_{\epsilon_i}^{N-\epsilon_i}(x)dx =0. $, on en d\'eduit gr\'ace au proc\'ed\'e d'it\'eration de Moser que $ \sup_{B(y,\delta_y/2)} u_{\epsilon_i} \to 0 $ et donc $ ||u_{\epsilon_i}||_{L^{\infty}(K)}\to 0 $. Ainsi, $ ||u_{\epsilon_i}||_{L^{\infty}(\bar \Omega)}\to 0 $, or ceci contredit le fait que $ ||u_{\epsilon_i}||_{L^N(\Omega)}\geq c>0 $, pour tout $ i $.

\bigskip

Pour la suite $ (u_{\epsilon_i}) $, on a un nombre fini de points de concentration de $ \Omega $, on consid\`ere alors, une sous-suite $ (u_{\epsilon_j}) $ de $ (u_{\epsilon_i}) $ qui a le maximum de points de concentration, soit $ x_1,x_2,\ldots,x_m $ tous les points de concentration de $ (u_{\epsilon_j}) $. Alors toute sous-suite de $ (u_{\epsilon_j}) $ a au plus $ m $ points de concentrations.

\bigskip

D'apr\`es le lemme 1( modifi\'e, $ \limsup \to \liminf $) de [He], on a:

$$ \forall \,\, k \in \{ 1,\ldots,m \} \,\, \forall \,\, \delta >0 \,\, \liminf_{\epsilon_j \to 0} \int_{B(x_k,\delta)} u_{\epsilon_j}^{N-\epsilon_j}(x)dx \geq \dfrac{\omega_n}{2^n}. $$

et,

$$ \forall \,\, x\in \Omega-\{x_1,\ldots,x_m\} \,\,\exists \,\, \delta_y>0 \,\, \liminf_{\epsilon_j \to 0} \int_{B(y,\delta_y)} u_{\epsilon_j}^{N-\epsilon_j}(x)dx =0. $$

Soit $ K $ un compact de $\bar \Omega-\{x_1,\ldots,x_m \}$, on a, $ \sup_K u_{\epsilon_j}=u_{\epsilon_j}(y_j) $ avec, par compacit\'e, $ y_j \to y \in K $. En utlisant le m\^eme raisonnement que [He] dans le lemme 4 et le lemme 8 et les estimations elliptiques (th\'eor\`eme 9.11 de [GT]), puis les estimations de Schauder on en d\'eduit qu'il existe une sous-suite $ u_{\epsilon_{j_r}} $  qui tend vers 0 dans $ {\cal C}^2(K) $. De plus cette sous-suite v\'erifie:

$$ \forall \,\, k\in \{1,\ldots,m\} \,\, \forall \delta >0 \,\, \liminf_{\epsilon_{j_r} \to 0} \int_{B(x_k,\delta)} u_{\epsilon_{j_r}}^{N-\epsilon_{j_r}}(x)dx \geq \dfrac{\omega_n}{2^n}. $$

Ce qui revient \`a dire que la nouvelle sous-suite a au moins $ m $ points de concentrations, d'apr\'es la d\'efinition de $ m $, cette sous-suite a exactement $ m $ points concentrations $x_1,\ldots,x_m $. On obtient:

$$ \forall \,\, y\in \Omega-\{x_1,\ldots,x_m\} \,\,\exists \,\, \delta_y>0 \,\, \liminf_{\epsilon_{j_r} \to 0} \int_{B(y,\delta_y)} u_{\epsilon_{j_r}}^{N-\epsilon_{j_r}}(x)dx =0. $$

En consid\'erant une suite exhaustive de compacts $ (K_n) $ de $ \bar \Omega-\{ x_1,\ldots, x_m \} $, $ K_n \subset \dot K_{n+1} $, $ \cup_n K_n =\bar \Omega-\{x_1,\ldots,x_m\} $ et en utilisant le proc\'ed\'e diagonal, on en d\'eduit qu'il existe une sous-suite not\'ee $ u_{\epsilon_j} $ qui a exactement $ m $ points de concentrations $ x_1,\ldots, x_m $ et qui converge vers 0 dans $ {\cal C}_{loc}^2(\bar \Omega-\{x_1,\ldots,x_m\}) $.
\bigskip

Soit $ x_1 $ un point de concentration, alors, il existe une suite $ j_n \to +\infty $ et $ x_{j_n,1} \to x_1 $ avec $ u_{\epsilon_{j_n}}(x_{j_n,k}) \to+\infty $. En appliquant le m\^eme proc\'ed\'e pour $ x_2 $ puis pour $ x_k $, $ k $ allant de 3 \`a $ m $, on extrait une sous-suite not\'ee encore $ (u_{\epsilon_j}) $ et des suites de points $ (x_{j,k}) $, $ k $ allant de 1 \`a $ m $ tels que, $ u_{\epsilon_j}(x_{j,k}) \to +\infty $ et $ x_{j,k} \to x_k $.

\bigskip

Ainsi, pour tout $ k \in \{1,\ldots,m \} $, il existe une suite $ y_{j,k} \to x_k $ et $ u_{\epsilon_j}(y_{j,k}) \to +\infty $.

\bigskip

On pouvait remarquer que,

$$ |B(x_k,\delta)| \sup_{B(x_k,\delta)} u_{\epsilon_i}^{N-\epsilon_i} \geq \int_{B(x_k,\delta)} u_{\epsilon_i}^{N-\epsilon_i}(x)dx \geq \dfrac{\omega_n}{2^{n+1}}\,\,\, {\rm pour }\,\, i\geq i_0. $$

puis, on fait tendre $ \delta $ vers z\'ero.

\bigskip
 
\underbar {\bf Point ii):}

\bigskip
 
Soit $ x_k $ un point de concentration, alors:

$$ \mu_k=\liminf_{\epsilon_j \to 0} \int_{B(x_k,\delta_0)} u_{\epsilon_j}^{N-\epsilon_j}(x)dx \geq \dfrac{\omega_n}{2^n} , $$

avec, $ \delta_0=\min_{\{ i \not= j, i,j=1,\ldots,m \} } \dfrac{ d(x_i,x_j)}{2} $.

\smallskip

Alors, pour $ 0 < \delta \leq \delta_0 $, on a:

$$ \mu_k=\liminf_{\epsilon_j \to 0} \int_{B(x_k,\delta_0)} u_{\epsilon_j}^{N-\epsilon_j}(x)dx = \liminf_{\epsilon_j \to 0} \int_{B(x_k,\delta)} u_{\epsilon_j}^{N-\epsilon_j}(x)dx . $$

En consid\'erant une sous-suite et en prenant $ \phi \in {\cal C}(\Omega) $, on obtient:

$$ \int_{\Omega} u_{\epsilon_j}^{N-\epsilon_j}(x)\phi(x) dx=\sum_{k=1}^m \int_{B(x_k,\delta)} u_{\epsilon_j}^{N-\epsilon_j}(x)\phi(x)dx + \int_{\Omega-[\cup_{k=1}^m B(x_k,\delta)]} u_{\epsilon_j}^{N-\epsilon_j}(x)\phi(x)dx, $$

En prenant $ \delta $ assez petit et en utilisant le point i) du Th\'eor\`eme, on en d\'eduit que:

$$ \lim_{\epsilon_j \to 0} \int_{\Omega} u_{\epsilon_j}^{N-\epsilon_j}(x)\phi(x) dx = \sum_{k=1}^m \mu_k \phi(x_k) ,\,\,\, \mu_k\geq \dfrac{\omega_n}{2^n}. $$

\underbar {\bf Point iii):}

\bigskip

Soit $ x_{\epsilon_j} \in \Omega $ tel que $ u_{\epsilon_j}(x_{\epsilon_j})=\sup_{\Omega} u_{\epsilon_j} $, apr\'es extraction d'une sous-suite, on peut supposer que $ x_{\epsilon_j} \to x_0 \in \Omega $. Il est clair que $ x_0 $ est un point de concentration et donc $ x_0 $ est l'un des $ x_k $, $ k\in \{1,\ldots,m\} $, on peut supposer que $ x_0=x_1 $.

\bigskip

Soit $ K $ un compact de $ \Omega-\{x_1,\ldots,x_m\} $, alors il existe $ \alpha_1>0,\ldots, \alpha_m>0 $ tels que, $ K \subset \tilde \Omega = \Omega_1-[\cup_{i=1}^m { \bar B(x_i,\alpha_i)}] $, o\`u $ \Omega_1 $ est un ouvert relativement compact de $ \Omega $. D'apr\'es le point $ i) $, $ u_{\epsilon_j} \to 0 $ uniform\'ement sur $ {\bar {\tilde \Omega}} $.

\bigskip

Consid\'erons l'op\'erateur $ L=-\Delta+n(n-2){u_{\epsilon_j}}^{N-2-\epsilon_j} $, alors, $ L u_{\epsilon_j}=0 $ et l'in\'egalit\'e de Harnack usuelle est v\'erfi\'ee pour cet op\'erateur (voir [GT]). Ainsi, pour tout $ \hat \Omega \subset \subset \tilde \Omega $, il existe $ c=c(\hat \Omega,\tilde \Omega,n)>0 $ telle que:

$$ \sup_{\hat \Omega} u_{\epsilon_j} \leq c\inf_{\hat \Omega} u_{\epsilon_j}  . $$

On prend $ \hat \Omega $ contenant  $ K $ et tel que son bord ext\'erieur (bord le plus proche de $ \partial \Omega $ ) soit le bord d'un ouvert $ {\cal O} $ contenant une boule de centre $ x_1 $ et de rayon $ \alpha >0 $. Comme $ u_{\epsilon_j} $ est sous-harmonique, on obtient:

$$ \inf_{\hat \Omega} u_{\epsilon_j} \leq \inf_{\partial {\cal O}} u_{\epsilon_j}=\inf_{{\cal O}} u_{\epsilon_j}. $$

Comme, $ {\bar B(x_1,\alpha)} \in {\cal O} $, en utilisant l'in\'egalit\'e de Harnack du type $ \sup \times \inf $, on a:

$$ \sup_{B(x_1,\alpha)} u_{\epsilon_j} \times \inf_{{\cal O}} u_{\epsilon_j} \leq c=c({\cal O},\alpha,n). $$

Or, $ \sup_{\Omega} u_{\epsilon_j}=\sup_{B(x_1,\alpha)} u_{\epsilon_j} $, en combinant ces in\'egalit\'es, on obtient:

$$ \sup_{\Omega} u_{\epsilon_j} \sup_K u_{\epsilon_j} \leq c=c(K,\Omega,n). $$

\underbar {\bf Point vi):}

\smallskip

{ \underbar {\bf M\'ethode 1}}: 

\smallskip

On sait (voir [H] pages 164-165) qu'il existe un voisinage du bord $ \omega $, un ouvert $ \Omega' \subset \subset \Omega $ et une constante positive $ c'=c'(n,\Omega) $ tels que:

$$ \sup_{\omega} u_{\epsilon_j} \leq c' \int_{\Omega'} u_{\epsilon_j} (x) dx. $$

On peut (quitte \`a grandir $ \Omega' $), supposer qu'il existe $ R>0 $ tel que, $ B(x_k,R) \subset \Omega',\,\, k\in \{1,\ldots,m\} $.

\bigskip

Comme $ u_{\epsilon_j} $ est sous-harmonique d'apr\`es l'in\'egalit\'e de Moser-Harnack, il existe $ 0 < R'<R $  et une constante positive $ \tilde c=\tilde c(n,\Omega,R') $ tels que:

$$ \int_{B(x_k,R')} u_{\epsilon_j}(x)dx \leq \tilde c \inf_{B(x_k,R')} u_{\epsilon_j}=\tilde c \inf_{\partial B(x_k,R')} u_{\epsilon_j}, \,\, k \in \{1,\ldots,m\}. $$ 

D'apr\`es le point iii) , on a:

$$ \sup_{\Omega} u_{\epsilon_j} \times \sup_{\partial B(x_k,R')} u_{\epsilon_j} \leq \tilde c_k . $$

Ce qui donne:

$$ \sup_{\Omega} u_{\epsilon_j} \int_{B(x_k,R')} u_{\epsilon_j}(x)dx \leq \bar c(R',\Omega, k,n). $$

On se place maintenant sur $ \Omega'-\{ \cup_{k=1}^m B(x_k,R'/2)\} $. Sur cet ouvert on a d'apr\`es le point i), $  u_{\epsilon_j} \leq 1 $ \`a partir d'un certain rang, on peut appliquer l'in\'egalit\'e de Harnack usuelle \`a la fonction $ u_{\epsilon_j} $ solution de $ Lu_{\epsilon_j}=0 $ o\`u $ L=-\Delta +n(n-2) u_{\epsilon_j}^{N-2-\epsilon_j} $. Il existe $ \hat   c=\hat c(R',\Omega,n) $ telle que:

$$ \sup_{\Omega'-\{ \cup_{k=1}^m B(x_k,R')\} } u_{\epsilon_j} \leq \hat c \inf_{ \Omega'- \{ \cup_{k=1}^m B(x_k,R')\} } u_{\epsilon_j} = \hat c \inf_{\partial (\Omega'- \{ \cup_{k=1}^m B(x_k,R')\})} u_{\epsilon_j} \leq \hat c \inf_{ \partial \Omega' } u_{\epsilon_j}= \hat c \inf_{ \Omega'} u_{\epsilon_j} . $$ 

En utilisant l'in\'egalit\'e de Harnack du type $ \sup \times \inf $, on obtient:

$$ \sup_{\Omega} u_{\epsilon_j} \times \inf_{\Omega'} u_{\epsilon_j} =\sup_{B(x_1,R')} u_{\epsilon_j} \times \inf_{\Omega'} u_{\epsilon_j} \leq c(\Omega', \Omega, n). $$

Finalement, on obtient:

$$\sup_{\Omega} u_{\epsilon_j} \times \int_{\Omega'-\{ \cup_{k=1}^m B(x_k,R')\}} u_{\epsilon_j}(x)dx \leq  \sup_{\Omega} u_{\epsilon_j} \times \sup_{\Omega'-\{ \cup_{k=1}^m B(x_k,R')\} } u_{\epsilon_j} \leq c(\Omega,R',\Omega,n). $$

En combinant, toutes ces in\'egalit\'es, on obtient,

$$ \sup_{\Omega} u_{\epsilon_j} \times \sup_{\omega} u_{\epsilon_j}  \leq c' \sup_{\Omega}u_{\epsilon_j} \int_{\Omega'} u_{\epsilon_j} (x)dx \leq c(\omega, \Omega, n). $$

{\underbar {\bf M\'ethode 2:}}

\smallskip

D'apr\'es l'in\'egalit\'e d'Alexandrov-Bakelman-Pucci (voir [Au]), on a:

$$ \sup_{\omega} u_{\epsilon_j} \leq \sup_{\partial \omega} u_{\epsilon_j}+C||u_{\epsilon_j}^{N-1-\epsilon_j}||_{L^n(\omega)}. $$

Donc,

$$ \sup_{\omega} u_{\epsilon_j} \leq \sup_{\partial \omega} u_{\epsilon_j}+C||u_{\epsilon_j}^{N-1-\epsilon_j}||_{L^n(\omega)} \leq \sup_{\partial \omega} u_{\epsilon_j}+C(\sup_{\omega} u_{\epsilon_j})^{N-1-\epsilon_j} .$$

D'apr\'es le point i), $ \sup_{\omega} u_{\epsilon_j} \to 0 $, d'o\`u:

$$ (1-k)\sup_{\omega} u_{\epsilon_j} \leq \sup_{\partial \omega} u_{\epsilon_j} \,\, {\rm avec} \,\, 0 < k < 1 . $$

En utilisant iii), on conclut que:

$$ \sup_{\Omega} u_{\epsilon_j} \times \sup_{\omega} u_{\epsilon_j} \leq \dfrac{1}{1-k} \sup_{\Omega} u_{\epsilon_j} \sup_{\partial \omega} u_{\epsilon_j} \leq c(\omega,\Omega,n). $$

\underbar {\bf Point v):}

\bigskip
 
D'apr\`es la formule de Pohozaev (voir [P] ou [H]), on a:

$$ \epsilon_j \int_{\Omega} u_{\epsilon_j}^{N-1-\epsilon_j}(x)dx=c_n \int_{\partial \Omega} <(x-y)|\nu(x)>  \left [\partial_{\nu} u_{\epsilon_j}(\sigma) \right ]^2 d\sigma. $$

D'o\`u,

$$ \epsilon_j \left ( \sup_{\Omega} u_{\epsilon_j} \right )^2 \int_{\Omega} u_{\epsilon_j}^{N-1-\epsilon_j}(x)dx=c_n \int_{\partial \Omega} <x|\nu(x)>  \left [\partial_{\nu} [ ( \sup_{\Omega} u_{\epsilon_j}) u_{\epsilon_j}](\sigma) \right ]^2 d\sigma. $$

Posons, $ g_{\epsilon_j}=(\sup_{\Omega} u_{\epsilon_j}) u_{\epsilon_j} $ et $ f_{\epsilon_j}=n(n-2)(\sup_{\Omega} u_{\epsilon_j}) u_{\epsilon_j}^{N-1-\epsilon_j} $. Alors, on a:

$$ \Delta g_{\epsilon_j}=f_{\epsilon_j} \,\, {\rm dans} \,\, \Omega \,\, {\rm et} \,\, g_{\epsilon_j}=0 \,\, {\rm sur} \,\, \partial \Omega. $$

D'apr\`es le point iv) , $ ||g_{\epsilon_j}||_{L^{\infty}(\omega)} \leq c(\omega,\Omega,n) $, avec $ \omega $ un voisinage du bord $ \partial \Omega $.

\bigskip

D'apr\`es la formule de repr\'esentation de la fonction Green, on a:

$$ g_{\epsilon_j}(x)=\int_{\Omega} G(x,y)f_{\epsilon_j}(y)dy \geq \int_{\Omega-\omega} G(x,y)f_{\epsilon_j}(y) dy. $$

D'apr\`es le point  iii) , on a, $ g_{\epsilon_j}(x) \leq c(x,\Omega,n) $ pour $ x \in \Omega-\{x_1,\ldots,x_m\} $ et d'apr\'es le point $ iv) $ et $ i) $, $ ||f_{\epsilon_j}||_{L^{\infty}(\omega)} \leq c(\omega,\Omega,n) $. 

\bigskip

D'apr\'es le principe du maximum, $ G(x,y)\geq c(x,\Omega-\omega,\Omega,n) >0 $ pour $ x\in \Omega $.

\bigskip

Finalement,

$$ \int_{\Omega } f_{\epsilon_j}(y)dy \leq c\,\, \forall \,\, j.$$

On veut prouver qu'il existe une constante $ c=c(\Omega,n)>0 $ telle que pour toute solution positive $ u_i $
 de,
 
 $$ \Delta u_i={u_i}^{N-1-\epsilon_i},\,\, {\rm dans }\,\, u_i=0 \,\, {\rm sur }\,\, \partial \Omega, $$
 
 on a,
 
 $$ \int_{\Omega} f_i \leq c, $$

 On sait qu'il existe un voisinage  $ \omega $ du bord tel que $ \sup_{\omega} u_i \leq c(\omega,\Omega,n), $ on consid\`ere l'op\'erateur $ L=-\Delta+n(n-2){u_i}^{N-2-\epsilon_i} $ dans l'ouvert $ \omega $, pour cet op\'erateur on peut appliquer l'in\'egalit\'e de Harnack usuelle ( voir [GT], chapitre 8). On obtient,
 
 $$ \sup_{\tilde \omega} u_i \leq \hat c \inf_{\tilde \omega} u_i, $$
 
 pour tout ouvert $ \tilde \omega $ relativement compact dans $ \omega $ et $ \hat c =c(\tilde \omega,\omega,n) $.
 
 \bigskip
 
 On sait aussi, qu'il existe $ \delta=\delta(\Omega,n)>0 $ tel que $ d(x_i,\partial \Omega)\geq \delta $, avec  $ u_i(x_i)=\sup_{\Omega} u_i $. On prend alors, $ \tilde \omega $ tel que son bord ext\'erieur soit $ \partial \tilde \Omega $, et $ \tilde \Omega $ contient strictement les $ x_i $. On a $ \omega $ d\'epend de $ \Omega $ et de $ n $, de m\^eme $ \tilde \omega \subset \subset \omega $, est choisi \`a partir de $ \omega $ et donc d\'epend ( il le doit) de $ \Omega $ et de $ n $.
 
 $$ \sup_{\Omega} u_i \times \sup_{\tilde \omega } u_i\leq \hat c  \sup_{\Omega} u_i  \inf_{\tilde \omega} u_i \leq \hat c \sup_{\Omega} u_i \inf_{\partial_{exterieur}\tilde \omega} u_i=\hat c u_i(x_i)\times \inf_{\partial \tilde \Omega} u_i=\hat c u_i(x_i)\times \inf_{\tilde \Omega} u_i, $$
 
 car, $ u_i $ est sous-harmonique et l'inf sur l'ouvert est atteint sur le bord (principe du maximum).
 
 \bigskip
 
 On utilise l'in\'egalit\'e de Harnack du type $ \sup \times \inf $ pour obtenir,
 
 $$ \sup_{\Omega} u_i \times \inf_{\tilde \Omega } u_i\leq c(\tilde \Omega,\tilde \omega,n), $$

Donc,

$$ g_i(x)\leq c(\tilde \omega, \Omega,n)\,\,\, {\rm sur}\,\, \tilde \omega. \qquad (*) $$

D'apr\'es, l'in\'egalit\'e d'Alexandrov-Bekelman-Pucci,

$$  \sup_{\omega} u_i \leq \sup_{\partial \omega} u_i+C||u_{i}^{N-1-\epsilon_i}||_{L^n(\omega)} \leq \sup_{\partial \omega} u_i +\tilde C ||u_i||_{L^n (\omega)}, $$

En multipliant par $ \sup_{\Omega} u_i $ \`a droite et \`a gauche, on obtient,

$$ \sup_{\omega} g_i \leq \sup_{\partial \omega} g_i+C||g_i||_{L^n(\omega)} , $$

On peut \'ecrire, 

$$ || g_i||_{L^n(\omega)} \leq || g_i||_{L^n(\omega-\omega_{\epsilon})}+|| g_i||_{L^n(\omega_{\epsilon})}, $$

avec, $ \omega_{\epsilon} $ un ouvert de $ \omega $  qui est voisinage de $ \partial \Omega $,

$$ || g_i||_{L^n(\omega_{\epsilon})} \leq \sup_{\omega_{\epsilon}} g_i |\omega_{\epsilon}|^{1/n}\leq \sup_{\omega} g_i |\omega_{\epsilon}|^{1/n}, $$

ainsi,

$$ \sup_{\omega} g_i(1-\tilde C|\omega_{\epsilon}|^{1/n})\leq \sup_{\partial \omega} g_i+ \bar C \sup_{\omega-\omega_{\epsilon}} g_i, $$

En prenant $ \omega_{\epsilon} $ voisin de $ \partial \Omega $ de mesure tr\'es petite $|\omega_{\epsilon}|<\tilde C^{-n} $, puis en utilisant $ (*) $ avec $ \tilde \omega=\omega-\omega_{\epsilon} $  puis $ \tilde \omega=\partial \omega $, on obtient,

$$ || g_i||_{L^{\infty}}(\omega)\leq c(\Omega,n). $$

Gr\^ace \`a la fonction de Green du laplacien,

$$ g_i(x)=\int_{\Omega} G(x,y)f_i(y)dy \geq \int_{\Omega- \omega} G(x,y)f_i(y)dy, $$

Par le principe du maximum, $ G(x,y)\geq c(\partial \omega-\partial \Omega,\Omega-\omega, \Omega, n)) $ pour $ x\in \partial \omega-\partial \Omega $ et $ y\in \Omega-\omega $.

\bigskip

Donc,

$$ \int_{\Omega- \omega} f_i\leq c(\Omega,n), $$

or, $ f_i=n(n-2) g_i\times { u_i}^{N-2-\epsilon_i}\leq c(\Omega,n) $ sur $ \omega $, d'o\`u le r\'esultat.
 
\bigskip

En utilisant, le lemme 8 de [He] ou le lemme 2 de [H],il existe une constante positive $ c=c(\omega,\Omega,n) $ telle que pour $ \beta \in ]0,1] $, pour tout voisinage $ \omega' \subset \subset \omega $ du bord $ \partial \Omega $, on a:

$$ ||g_{\epsilon_j}||_{W^{1,q}(\Omega)}+||\nabla g_{\epsilon_j}||_{{\cal C}^{0,\beta}(\omega')} \leq c, \,\, 1\leq q<\dfrac{n}{n-2}. $$

En utilisant le th\'eor\`eme d'Ascoli, on peut extraire de $ g_{\epsilon_j} $ une sous-suite qu'on note encore $ g_{\epsilon_j} $ qui converge uniform\'ement dans $ {\cal C}^1(\omega') $ vers $ g\geq 0 $. D'apr\'es le principe du maximum de Hopf, $ \sup_{\partial \Omega} |\partial_{\nu} g| >0 $, en utilisant le point 1) et l'identit\'e de Pohozaev (voir [H]), on conclut qu'il existe deux constantes positives $ \beta_1 $ et $ \beta_2 $ telles que::

$$ \beta_1 \leq \epsilon_j \left ( \sup_{\Omega} u_{\epsilon_j} \right )^2 \leq \beta_2. $$

Plus pr\'ecis\'ement,

$$ \epsilon_j \left ( \sup_{\Omega} u_{\epsilon_j} \right )^2 \to \dfrac{ c_n \int_{\partial \Omega} <x|\nu(x)> [ \partial_{\nu} g(\sigma)]^2 d\sigma }{\sum_{k=1}^m \mu_k}. $$

%La fonction $ g $ v\'erifie au sens des distributions:

%$$ \Delta g=\sum_{k=1}^m \gamma_k \delta_{x_k},\,\, {\rm dans} \,\, \Omega \,\, {\rm et}\,\, g=0 \,\, {\rm sur} \,\, \partial \Omega, $$ 

%avec $ \gamma_k=\liminf_{\epsilon_j \to 0}\int_{B(x_k,\delta_0)} f_{\epsilon_j}(x)dx $  et $ \delta_0=\min_{i \not = j, i,j \in \{1,\ldots,m\}} \dfrac{d(x_i,x_j)}{2}. $ ($ f_{\epsilon_j} \to 0 $ uniform\'ement sur tout compact et donc $ \liminf_{\epsilon_j \to 0}\int_{B(x_k,\delta)} f_{\epsilon_j}(x)dx=\liminf_{\epsilon_j \to 0}\int_{B(x_k,\delta_0)} f_{\epsilon_j}(x)dx $  pour $ 0 < \delta <\leq \delta_0 $).

\bigskip

Sur un voisinage $ \omega $ de $ \partial \Omega $, en passant aux sous-suites, on peut supposer que $ g_{\epsilon_j} $ converge uniform\'ement vers $ g_{\omega} $ dans  $ {\cal C}^2(\omega) $.

\bigskip

\underbar {\bf Point vi):}

\bigskip

Consid\'erons une suite exhaustive de compacts $ (K_n) $ de $ \bar \Omega-\{x_1,\ldots,x_m\} $. Sur tout compact $ K $ de  $ \bar \Omega-\{x_1,\ldots,x_m\} $ le point iii) donne $ || g_{\epsilon_j}||_{L^{\infty}(K)} \leq c(K,\Omega,n) $. En utilisant les estimations elliptiques et le lemme 2 de [H] ou le lemme 8 de [He] et le th\'eor\`eme d'Ascoli, on peut extraire de $ (g_{\epsilon_j}) $ une sous-suite qui converge dans $ {\cal C}^1(K) $ vers une fonction $ g_K $. Comme $ \nabla f_{\epsilon_j}=n(n-2)(N-1-\epsilon_j) u_{\epsilon_j}^{N-2-\epsilon_j} \nabla g_{\epsilon_j} $, en appliquant ce proc\'ed\'e \'a la suite de compacts $ (K_n) $ et en utilisant le proc\'ed\'e diagonal et l'unicit\'e des limites, on construit de proche en proche une fonction $ g \in {\cal C}^2 (\bar \Omega-\{x_1,\ldots,x_m\}) $ et une sous-suite not\'ee encore $ g_{\epsilon_j} $ qui converge vers $ g $ dans $ {\cal C}_{loc}^2(\bar \Omega-\{x_1,\ldots,x_m\}) $. En prolonge $ g $ et $ \nabla g $ en $ \tilde g $ et $ \tilde \nabla g $ sur  $ \bar \Omega $ par $ 0 $ (par exemple) sur $ \{x_1,\ldots,x_m\} $. $ \tilde g $ reste une fonction mesurable sur $ \bar \Omega $ et comme $ || g_{\epsilon_j}||_{W^{1,q}(\Omega)} \leq c(n,\Omega) $ avec $ 1 < q <\dfrac{n}{n-1} $. Gr\^ace \`a l'injection de Sobolev $ W^{1,q} $ dans $ L^r $ $ 1\leq r\leq q^* $ et l'injection compacte de Khondrakov pour $ 1\leq r<q^* $, on a $ g_{\epsilon_j} \to \bar g $ dans $ L^r $ et $ \bar g=\tilde g $ presque partout( on compare $ g $ et $ \bar g $ sur $ K_n $, puis on passe \'a la limite en $ n$). Donc, $ \tilde g \in L^r(\Omega) $ avec $ r>1 $ et est solution ( au sens des distributions) de:

% on a $ \tilde g \in W_0^{1,q}(\Omega) $. En effet les prolongements $ \tilde g $ et $ \tilde \nabla g $ sont mesurables et comme $ ||g_{\epsilon_j}||_{W^{1,q}(K_r)} \leq ||g_{\epsilon_j}||_{W^{1,q}(\Omega)}\leq c(q,n,\Omega) $, on passe a la limite en $ j $ on obtient, $ ||\tilde g||_{L^q(K_r)}+||\nabla g||_{L^q(K_r)} \leq c(n,q,\Omega) $ puis on utilise le th\'eor\`eme de Beppo-Levi \`a la suite $ g_r=1_{K_r}g $. On prouve de la m\^eme mani\`ere que $ \tilde g \in W_0^{1,q} $, pour $ \phi \in {\cal D}(\Omega) $, on a :

%$$ \int_{\Omega} \partial_k g_{\epsilon_j}\phi=-\int_{\Omega} g_{\epsilon_j} \partial_k \phi , $$

%$$ \left | \int_{K_r} \partial_k g_{\epsilon_j}\phi+ g_{\epsilon_j} \partial_k \phi \right |\leq \int_{\Omega-K_r} |\partial_k g_{\epsilon_j}||\phi|+\int_{\Omega-K_r} |g_{\epsilon_j}||\phi|, $$

%En utilisant l\'in\'egalit\'e de H\"older:

%$$ \left | \int_{K_r} \partial_k g_{\epsilon_j}\phi+ g_{\epsilon_j} \partial_k \phi \right |\leq ||\phi||_{{\cal C}^1(\Omega)}[||g_{\epsilon_j}||_{W^{1,q}(\Omega)}] [mes(\Omega-K_r)]^{1/q'}. $$

%En passant \`a la limite en $ j $ puis en $ r $ (avec Beppo-Levi), on obtient:

%$$ \int_{\Omega-\{x_1,\ldots,x_m\}} \partial_k g\phi=-\int_{\Omega-\{x_1,\ldots,x_m\} } g\partial_k \phi. $$ 

%Apr\'es prolongement, on obtient $ \tilde g \in W_0^{1,q}(\Omega) $ et v\'erifie:

$$ \int_{\Omega} \tilde g \Delta \phi=\sum_{k=1}^m \gamma_k\delta_{x_k}\,\, \forall \,\, \phi \in \, {\cal C}^2(\bar \Omega)\cap \dot H_1^2(\Omega), $$

avec, $ \gamma_k=\liminf_{\epsilon_j \to 0}\int_{B(x_k,\delta_0)} f_{\epsilon_j}(x)dx $  et $ \delta_0=\min_{i \not = j, i,j \in \{1,\ldots,m\}} \dfrac{d(x_i,x_j)}{2}$. En utilisant les points i), iii) et iv),la suite $ (f_{\epsilon_j}) $ converge uniform\'ement vers 0 sur tout compact de $ \bar \Omega-\{x_1,\ldots,x_m\} $.

\bigskip

Or, la fonction G  de Green du laplacien avec condition de Dirichlet, v\'erifie:

$$ \int_{\Omega} G(x_k,x) \Delta \phi=\delta_{x_k}\,\, \forall \,\, \phi \in \, {\cal C}^2(\bar \Omega)\cup \dot H_1^2(\Omega). $$

Donc,

$$ \int_{\Omega} [\tilde g-\sum_{k=1}^m \gamma_k G(x_k,x)] \Delta \phi=0\,\, \forall \,\, \phi \in \, {\cal C}^2(\bar \Omega)\cap \dot H_1^2(\Omega) $$

et,

$$ \tilde g-\sum_{k=1}^m \gamma_k G(x_k,x)=0 \,\, {\rm sur } \,\, \partial \Omega. $$

En utilisant le th\'eor\`eme de r\'egularit\'e d'Agmon (voir th\'eor\`eme 8.2  pp 444 de [Ag]ou [Au]), on conclut que, $ u=\tilde g-\sum_{k=1}^m \gamma_k G(x_k,x) $ presque partout avec $ u\in {\cal C}^{\infty}(\Omega) $. 

\smallskip

Comme $ \tilde g-\sum_{k=1}^m \gamma_k G(x_k,.) $ est $ {\cal C}^1(\bar \Omega-\{x_1,\ldots,x_m\}) $, on en d\'eduit que $ u=g-\sum_{k=1}^m \gamma_k G(x_k,.) $ partout sur  $ \Omega-\{x_1,\ldots,x_m\}$. 

\smallskip

Soit alors, $ w $ la fonction qui vaut $ u $ dans $ B(x_k,\epsilon), k=1...m $ et $ g-\sum_{k=1}^m \gamma_k G(x_k,.) $ partout ailleurs sur $ \bar \Omega $, $ w $ est ${\cal C}^2(\bar \Omega) $ \'egale \'a $ u-\sum_{k=1}^m \gamma_k G(x_k,.) $ presque partout sur $ \Omega $ et donc, $ \int_{\Omega} w \Delta \phi=0 $ au sens des distributions et $ w\equiv 0 $ sur le bord $ \partial \Omega $. Finalement $ w\equiv 0 $ sur $ \Omega $; d'o\'u $ \tilde g\equiv \sum_{k=1}\gamma_k G(x_k,.) $ sur $ \bar \Omega -\{x_1,x_2,\ldots,x_m\} $. Ainsi, on a:

$$ \sup_{\Omega} u_{\epsilon_j} u_{\epsilon_j} \to \sum_{k=1}^m \gamma_k G(x_k,.) \,\, {\rm dans}\,\, {\cal C}_{loc}^2(\bar \Omega-\{x_1,\ldots,x_m\}). $$

On a:

$$\gamma_k=\liminf_{\epsilon_j \to 0} \int_{B(x_k,\delta_0) }f_{\epsilon_j} \,\, {\rm et} \,\,\liminf_{\epsilon_j \to 0} \int_{B(x_k,\delta_0)} u_{\epsilon_j}^{N-\epsilon_j}=\mu_k \geq \dfrac{\omega_n}{2^n},$$

donc,

$$ \gamma_k=n(n-2)\liminf_{\epsilon_j \to 0}   \int_{B(x_k,\delta_0)} \left [ (\sup_{\Omega} u_{\epsilon_j}) u_{\epsilon_j}^{N-1-\epsilon_j} \right ]  \geq n(n-2)\mu_k. $$

\newpage

{\underbar {\bf Preuve du Th\'eor\`eme 4:}}

\smallskip

Les conditions $ b,c>0 $, et $  c\leq \dfrac{ b^2 }{4} $, nous permettent d'avoir l'existence d'une fonction de Green $ G $ pour l'op\'erateur $ \Delta ^2 + b \Delta + c $, telle que,

$$ \dfrac{C'(M,g)}{d_g(x,y)^{n-4}}\geq G(x,y)\geq \dfrac{C(M,g)}{d_g(x,y)^{n-4}},$$

avec, $ C(M,g),C'(M,g)>0 $.

\bigskip

On \'ecrit alors,

$$ \min_M u_i=u_i(x_i)=\int_M G(x,y)V_i(y){u_i(y)}^{(n+4)/(n-4)}dV_g(y), $$

et,

$$ \min_M u_i \geq \bar C \int_M V_i(y){u_i(y)}^{(n+4)/(n-4)}dV_g(y), $$

d'o\`u,

$$ \sup_M u_i \times \inf_M u_i \geq \int_M V_i {u_i}^{2n/(n-4)}. $$

On multiplie l'\'equation $ (E) $ par $ u_i $ puis on int\`egre par parties, on obtient,

$$ \int_M V_i {u_i}^{2n/(n-2)} = \int_M [ (\Delta u_i)^2+ b|\nabla u_i|^2+c{u_i}^2 ] \geq K ||u_i||_{H_2^2(M)}^2, $$

On utilise l'injection de Sobolev $ H_2^2(M) $ dans $ L^{2n/(n-4)}(M) $ et le fait que $ 0 \leq V_i(x) \leq A $ sur $ M $ pour obtenir,

$$ A||u_i||_{L^{2n/(n-4)}}^{2n/(n-4)} \geq K'||u_i||_{L^{2n/(n-4)}}^{2}, $$

donc,

$$ ||u_i||_{L^{2n/(n-4)}} \geq K''> 0,\,\,\, \forall \,\, i $$

ainsi,

$$ \sup_M u_i \times \inf_M u_i \geq \int_M V_i {u_i}^{2n/(n-2)} \geq \tilde K >0 \,\, \forall \,\, i. $$

{\underbar {\bf Preuve du Th\'eor\`eme 5:}}

Supposons par l'absurde que:

$$ \sup_{\Omega} u_i \times \inf_K u_i \to 0. $$

Alors, pour $ \delta >0 $ assez petit, on a:

$$ \sup_{\Omega} u_i \times \inf_{\{x, d(x,\partial \Omega) \geq \delta \}} u_i \to 0 . $$

D'apr\`es la preuve de la proposition 2.4 de [C-G], on a pour $ \delta>0 $ assez petit,

$$ \sup_{ \{x, d(x,\partial \Omega) \leq \delta \}} u_i \leq M=M(n,\Omega). $$ 

On a,

$$ u_i(x)=\int_{\Omega} G(x,y){u_i(y)}^{p-\epsilon_i} dy, $$

Soient, $ K' $ un autre compact de $ \Omega $, en utilisant le principe du maximum, on obtient:

$$ \exists \,\,\, c_1=c_1(K, K', n, \Omega)>0,\,\,\,{\rm tel \, que}\,\, G(x,y)\geq c \,\,\, \forall \,\, x\in K,\,\, y\in K' , $$

donc,

$$ \inf_K u_i=u_i(x_i)\geq c_1 \int_{K'} {u_i(y)}^{p-\epsilon_i}dy, $$

En prenant, $ K=K_{\delta}=\{x,d(x,\partial \Omega)\geq \delta \} $, il existe $ c_2=c_2(\delta, n, K, \Omega)>0 $ telle que:

$$ \inf_K u_i \geq c_2 \int_{K_{\delta}} {u_i(y)}^{p-\epsilon_i} dy, $$

d'o\`u,

$$ \sup_{\Omega} u_i \times \inf_K u_i \geq c_2 \int_{K_{\delta}} {u_i(y)}^{p+1-\epsilon_i} dy . $$

On en d\'eduit que:

$$ ||u_i||_{p+1-\epsilon_i}^{p+1-\epsilon_i}=\int_{\{x, d(x,\partial \Omega) \leq \delta \}} {u_i}^{p+1-\epsilon_i}+\int_{\{x, d(x,\partial \Omega) \geq \delta \}} {u_i}^{p+1-\epsilon_i}, $$

Ce qui donne, 

$$ ||u_i||_{p+1-\epsilon_i}^{p+1-\epsilon_i} \leq \sup_{\Omega} u_i\times \inf_{\{x, d(x,\partial \Omega) \geq \delta \}} u_i+mes(\{x, d(x,\partial \Omega) \leq \delta \})M^{p+1-\epsilon_i}, $$

en faisant tendre, $ i $ vers l'infini et en prenant $ \delta $ assez petit, on conclut que,

$$ ||u_i||_{p+1-\epsilon_i} \to 0. $$

Or, d'apr\`es l'injection de Sobolev (voir [V]), $ H_2^2(\Omega)\cap H_0^1(\Omega) $ dans $ L^{p+1}(\Omega) $, en multipliant $ (E) $ par $ u_{\epsilon} $ , en int\`egrant par parties et en utilisant l'in\'egalit\'e de H\''older, on obtient,

$$ K_1||u_i||_{p+1-\epsilon_i}^2\leq K_2 ||u_i||_{p+1}^2 \leq \int_{\Omega} (\Delta u_i)^2=\int_{\Omega} {u_i}^{p+1-\epsilon_i}=||u_i||_{p+1-\epsilon_i}^{p+1-\epsilon_i}. $$

De $ 0 < \epsilon_i \leq \dfrac{4}{n-4} $ et de l'in\'egalit\'e pr\'ec\'edente, on a la contradiction suivante,

$$ ||u_i||_{p+1-\epsilon_i} \geq K_3>0, $$

\bigskip

{\underbar {\bf Remerciements}}

\bigskip

Ce travail \`a \'et\'e fait pendant le s\'ejour de l'auteur en Gr\`ece. L'auteur tiens \`a remercier le D\'epartement de Math\'ematiques de l'Universit\'e de Patras, surtout le Professeur Athanase Cotsiolis et la Fondation IKY pour leur acceuil.

\bigskip

{\underbar {\bf Bibliographie:}}

\bigskip

[Ag] S. Agmon. The $ L_p $ Approach to the Dirichlet Problem. Ann. Scuola Norm. Sup. Pisa 13, (1958) 405-448.

\smallskip

[Au] T. Aubin. Some nonlinear Problems in Riemannian Geometry. Springer Verlag 1998.

\smallskip

[Au, D, He] T.Aubin, O. Druet, E. Hebey, Best constants in Sobolev inequalities for compact manifolds of nonpositive curvature. C.R. Acad. Sci. Paris S\'er. I Math, 326 (1998), no 9 1117-1121.

\smallskip

[B 1] S.S Bahoura. Majorations du type $ \sup u \times \inf u \leq c $ pour l'\'equation de la courbure scalaire prescrite sur un ouvert de $ {\mathbb R}^n, n\geq 3 $. J.Math.Pures Appl.(9) 83 (2004), no.9, 1109-1150.

\smallskip

[B-2] S.S Bahoura. Estimations du type $ \sup  \times \inf $ sur une vari\'et\'e compacte (\`a paraitre).

\smallskip

[B 3] S.S Bahoura. In\'egalit\'es de Harnack pour les op\'erateurs d'ordre 2 et 4 et ph\'enom\`ene de concentration. C.R.A.S.

\smallskip

[BP] H. Br\'ezis, L. A. Peletier, Asymptotic for elliptic equations involving critical growth. Partial eifferential equations and calculus of variation, vol.I, 149-192, Progr.Nonlinear Differential Equations Appl.,1, Birkhauser Boston, Boston,MA,1989.

\smallskip

[C 1] D. Caraffa. Etude des probl\`emes elliptiques non lin\'eaires du quatrieme ordre avec exposants critiques sur les vari\'et\'es riemanninnes compactes. J. Math. Pures. Appl. (9) 83, (2004), no.1, 115-136.

\smallskip

[C 2] D. Caraffa. Equations elliptiques du quatri\`eme ordre avec exposants critiques sur les vari\'et\'es riemanniennes compactes. J. Math. Pures. Appl. (9) 80, (2001) no.9, 941-960.

\smallskip

[C-L] C-C Chen, C-S Lin. Estimates of the conformal scalar curvature equation via the mathod of moving planes. Comm. Pure Appl. Math. 37 (1997) 0971-1017.

\smallskip

[C-G] K-S. Chou, D. Geng. Asymptotics Of Positive Solutions For A Biharmonic Equation Involving Critical Exponent. Difftial and Integral Equations. Volume 13 (7-8) July-September 2000, pp. 921-940.

\smallskip

[DLN] D.G. De Figueiredo, P.L. Lions, R.D. Nussbaum, A priori Estimates and Existence of Positive Solutions of Semilinear Elliptic Equations, J. Math. Pures et Appl., vol 61, 1982, pp.41-63.

\smallskip

[D H R] O.Druet, E.Hebey, F.Robert, Blow-up theory for elliptic PDEs in Riemannian Geometry. Mathematical Notes, 45. Princeton University Press, Princeton, NJ, 2004.
\smallskip

[GNN] B. Gidas, W. Ni, L. Nirenberg, Symmetry and Related Propreties via the Maximum Principle, Comm. Math. Phys., vol 68, 1979, pp. 209-243.

\smallskip

[GT] D. Gilbarg, N.S. Trudinger. Elliptic Partial Differential Equations of Second order, Berlin Springer-Verlag, Second edition, Grundlehern Math. Wiss.,224, 1983.

\smallskip

[H] Z-C. Han, Assymptotic Approach to singular solutions for Nonlinear Elleptic Equations Involving Critical Sobolev Exponent. Ann. Inst. Henri Poincar\'e. Analyse Non-lin\'eaire. 8(1991) 159-174.

\smallskip

[He] E. Hebey. Asymptotics for some quasilinear elliptic equations. Diff and Int Eq. Volume 9, Number 1, (1996), pp. 71-88.

\smallskip

[He,V] E. Hebey, M. Vaugon, The best constant problem in the Soboblev embedding theorem for complete riemannian  Manifolds. Duke Math.J. 79 (1995), no 1, 235-279.

\smallskip

[L1] Y.Y Li. Prescribing Scalar Curvature on $ {\mathbb S}_n $ and related Problems. I. J. Differential Equations 120 (1995), no. 2, 319-410.

\smallskip

[L2] Y.Y Li. Prescribing Scalar Curvature on $ {\mathbb S}_n $ and related Problems. II. Comm. Pure. Appl. Math. 49(1996), no.6, 541-597.

\smallskip

[M] J. Moser. On Harnack's Theorem for Elliptic Differential Equations. Comm. Pure Appl Math. vol 15, 577-591 (1961).

\smallskip

[P] S. Pohozaev, Eigenfunctions of the equation $\Delta
 u+\lambda f(u)=0 $. Soviet. Math. Dokl., vol. 6 (1965), 1408-1411.

\smallskip

[T] N.S. Trudinger, Remarks Concerning The Conformal Deformation of Riemannian Structures On Compact Manifolds. Ann. Scuola Norm. Sup. Pisa. 22 (1968) 265-274.

\smallskip

[V] R.C.A.M. van der Vorst. Best constant for the embedding of the space $ H^2\cap H_0^1 $ into $ L^{2N/(N-4)}(\Omega) $. Diff. Int. Eq., 6 (1993), 259-276.

\smallskip

%[DLN] D. G. de Figueirdo, P-L. Lions, and R.D. Nussbaum, A priori estimates and existence of positive solutions of semilinear equations. J. Math. Pure Appl., 61 (1982), 41-63.
\smallskip

%[H] Z. Han. Asymptotic approach to singular solutions for nonlinear elliptic equations involving critical Sobolev exponent, Ann. Ins. Henri Poincar\'e Analyse non lin\'eaire. 8 (1991) 159-174.
\smallskip

%[V] R.C.A.M. van der Vorst. Best constant for the embedding of the space $ H^2\cap H_0^1 $ into $ L^{2N/(N-4)}(\Omega) $. Diff. Int. Eq., 6 (1993), 259-276.

\end{document}